\documentstyle[twoside,amssymb,12pt]{article}
\def\C{{\Bbb C}}

\def\Q{{\Bbb Q}}

\newtheorem{prop}{Proposition}%[section]
\newtheorem{dfn}[prop]{Definition}
\newtheorem{theo}[prop]{Theorem}

\newtheorem{lem}[prop]{Lemma}

\title{\sc Canonical abelianization of finite group actions} 

\author{{\sc Victor V. Batyrev} \\
\small  {\em Mathematisches Institut, Universit\"at T\"ubingen}   \\
\small  {\em Auf der Morgenstelle 10,  72076  T\"ubingen, Germany}  \\
\small  {\em e-mail: batyrev@bastau.mathematik.uni-tuebingen.de} \\
 }

\begin{document}

\date{}

\maketitle

\begin{abstract}
In these notes  we present a corrected version of the 
algorithm for abelianization 
of stabilizers of finite group actions on smooth manifolds 
using equivariant blow ups 
(see \cite{nonarch} \S 5) 
\end{abstract}

Let $G$ be a finite group, $V$ a smooth $n$-dimensional 
algebraic variety over ${\C}$  
having  a regular effective 
action of $G$. If $x \in V$ is an arbitrary 
 point,  then by $St_G(x)$ we denote the stabilizer of $x$ in $G$.
 For any subset   $M \in G$ we 
define $V^M := \{ x \in V\, : \, gx=x \;\; \forall g \in M\}$.

\begin{dfn} 
{\rm Let  $D = \sum_{i=1}^m d_i D_i \in {\rm Div}(V)^G \otimes \Q$ an 
effective  
$G$-invariant $\Q$-divisor on a $G$-manifold $V$.  A pair  
$(V,D)$ will be called  {\bf $G$-normal}  if the following 
conditions are satisfied: 

(i) $Supp\, D$ is a union of normal crossing divisors $D_1, \ldots, D_m$; 

(ii) for any irredicible component 
 $D_i$ of $D$ and for any point $x \in D_i$ 
 the divisor  $D_i$ is $St_G(x)$-invariant,
i.e., $g(D_i) = D_i$,  $\;\forall \, g \in 
St_G(x)$ (however, the $St_G(x)$-action on $D_i$ itself 
may be nontrivial).  
} 
\label{norm-p}
\end{dfn}

\begin{theo} \label{abel1}
Let $(V,D)$ be a $G$-normal pair. Then there exists a  
canonically determined sequence of smooth 
$G$-manifolds $V_0 := V, V_1, \ldots, V_k = V^{ab}$ 
and birational $G$-euivariant morphisms $\varphi_i\, : \, V_i \to V_{i-1}$ 
such that each $\varphi_i$ is a blow up of $V_{i-1}$ at a 
$G$-invariant submanifold and the birational 
morphism $\psi:= \varphi_k \circ \cdots \circ \varphi_1 \, 
:\, V^{ab} \rightarrow V$ 
has the following properties: 

{\rm (i)} Let $D^{ab} = (K_{V^{ab}} -  \psi^*K_{V})  +  \psi^{-1}D$, 
where   $\psi^{-1}D$ is a proper pullback of $D$ with respect to 
the birational morphism $\psi$, then $(V^{ab}, D^{ab})$
is a  $G$-normal pair; 

{\rm (ii)} for any point $x \in V^{ab}$  
the stabilizer $St_G(x)$ is an abelian subgroup in $G$. 
\label{canon-ab} 
\end{theo}

\noindent 
{\em Proof.} 
First of all we remark that from the algorithm presented 
below it follows   that 
for each $G$-equivariant blow up $\varphi_i\, : \, V_i \to V_{i-1}$
the divisor $D_i:= (K_{V_i} - \varphi^*K_{V_i}) + \varphi_i^{-1}D_{i-1}$
satisfies both conditions (i),(ii) in \ref{norm-p} if $(V_{i-1}, D_{i-1})$
is a $G$-normal pair. Since in  \ref{norm-p} only $Supp\, D$
is important, the property  \ref{canon-ab}(i) will hold authomatically.  
Thus, our main goal will be to get  \ref{canon-ab}(ii). 

Let $Z(V,G) \subset  V$ be the set of all points 
$x \in V$ such that  $St_G(x)$ is not abelian. If $Z(V,G)$ is empty, 
then we are 
done.  Assume that 
$Z(V,G) \neq \emptyset$. 
We set 
\[ s(V,G) := \max_{x \in Z(V,G)} |St_G(x)|. \]
Consider  a Zariski closed subset 
\[ Z_{\rm max}(V,G) := \{ x \in Z(V,G) \; :\;  |St_G(x)| = s(V,G) \} \subset 
Z(V,G). \]
We claim that the set  $Z_{\rm max}(V,G) \subset V$ is a smooth  $G$-invariant
subvariety of  codimension at least $2$. By definition, 
$Z_{\rm max}(V,G)$ is a  union of smooth subvarieties  
$$F(H):= \{ x \in V\; : \; gx =x\;\; \forall g \in H \},$$ 
where $H$ runs over all nonabelian 
subgroups of $G$ such that  $|H| =  s(V,G)$. 
This implies that $ Z_{\rm max}(V,G)$ is $G$-invariant.  
Since the $G$-action is effective and $dim\, F(H) =n-1$
is possible only for cyclic subgroups $H \subset G$. So we obtain  
$dim\, Z_{\rm max}(V,G) \leq n-2$.  
It remains to observe that  any two subvarieties $F(H_1), F(H_2) \subset 
V$ must either coincide, or have empty 
intersection. Indeed, if $x \in F(H_1) \cap  F(H_2)$, then 
$H_1, H_2 \subset St_G(x)$. From maximality of  $|H_1|$ and  $|H_2|$ 
it follows that  
$H_1 = H_2 =  St_G(x)$; i.e., $F(H_1) = F(H_2)$.

Denote by $\{ Z_1, \ldots, Z_k\}$ the set of all distinct 
connected components of  
$Z_{\rm max}(V,G)$. By  definition, 
of $Z_{\rm max}(V,G)$, we have $St_G(x) =  St_G(x')$ for 
any two points $x, x' \in Z_i$. Thus  we obtain  
a sequence of  maximal nonabelian stabilizers  
$H_1, \ldots, H_k \subset G$ : $H_i = St_G(x)$ $(x \in Z_i)$. We admit that 
$H_i$ may coincide with $H_j$ for $i \neq j$.  
Let us consider 
the linear representation $\rho_i\,:\, H_i \rightarrow GL(n, \C) = GL(TV_x)$, 
where $TV_x$ the tangent space of $V$ at $x$. It follows from the 
connectedness of $Z_i$ that $\rho_i$ does not depend on the choice of 
$x \in Z_i$. We write $\rho_i$ as a direct sum of irreducible 
representations 
$$ \rho_i = \chi^{(i)}_1 \oplus \cdots \oplus \chi^{(i)}_{r_i} \oplus 
\tau^{(i)}_1  \oplus \cdots \oplus \tau^{(i)}_{s_i},$$
where $\chi^{(i)}_1, \ldots, \chi^{(i)}_{r_i}$ are characters (representations 
of degree $1$) and 
each irreducible representation 
$\tau^{(i)}_j$ $(j =1, \ldots, s_i)$ has degree $\geq 2$. 
Denote by $\overline{H_i} \subset H$ the common 
kernel of the characters  $\chi^{(i)}_1, \ldots, \chi^{(i)}_{r_i}$ 
(we set $\overline{H_i} = H_i$ if $r_i = 0$).

Now we define $C(V,G)$ to be the union of those 
connected components $Z_i \subset Z_{\rm max}(V,G)$ such that $\overline{H_i}
 = H_i$ 
(the latter holds if all characters $\chi^{(i)}_1, \ldots, \chi^{(i)}_{r_i}$
are trivial or if $r_i = 0$). 
It follows from the above properties of  $Z_{\rm max}(V,G)$
that $C(V,G)$ is $G$-invariant, smooth and has codimension at 
least $2$. 

We distinguish the following cases:
\medskip 

{\sc Case 1}:  $C(V,G) \neq \emptyset$. Then we 
set $V_1$ to be the $G$-equivariant 
blow-up of $V_0$ 
with center $C(V,G)$. Denote by   
$\varphi_1 \, :\, V_1 \rightarrow V_0$ the corresponding projective birational 
$G$-morphism. It is obvious that the support of $D_1 := 
K_{V_1} - \varphi_1^*K_V + \varphi_1^{-1}D$ is a normal crossing divisor. If  
$E$ is a connected component of the 
$\varphi_1$-exceptional divisor and $x \in E$ is its arbitrary point, 
then $St_G(x) \subset St_G(\varphi_1(x))$. Since 
$\varphi(E) = Z_i$ is a connected component of $C(V,G)$, 
the stabilizer of each point $\varphi_1(x) \in Z_i$ 
is equal to $H_i$, and the divisor $E$ must be  $St_G(x)$-invariant. 
Hence, we conclude that $(V_1,D_1)$ is 
a $G$-normal pair. In the our situation  $St_G(x) \neq 
St_G(\varphi_1(x))$, because the $H_i$-action in the normal space to $Z_i$ at 
$\varphi_1(x)$ has no $1$-dimensional $H_i$-invariant subspaces 
($r_i = {\rm dim} \, Z_i$). So 
$|St_G(x)| < |St_G(\varphi(x))|$. From our 
construction it  follows that either $s(V_1, G) < s(V_0, G)$, or 
the number of connected components of  $Z_{\rm max}(V_1,G)$ is less than 
the  number of connected components of  $Z_{\rm max}(V_0,G)$. Using induction, 
we eventually come to the second case:  
\medskip

{\sc Case 2}:  $C(V,G) = \emptyset$. This means that $\overline{H_i} \neq H_i$ 
(in particular, $r_i > 0$) for all $i \in \{ 1, \ldots, k \}$. 
We remark that  all subgroups $\overline{H_i}$ are nontrivial, since all 
subgroups $H_i$ are nonabelian.  
For every  $i \in \{ 1, \ldots, k \}$ we define $W_i$ to be the 
connected component of $V^{\overline{H_i}}$ containing $Z_i$. 
\medskip

\begin{lem} 
One has  ${\rm dim}\, W_i = r_i$. In particular, the 
codimension of each $W_i$ is equal to 
$\sum_{l =1}^{s_i} {\rm deg}\, \tau_l \geq 2$.
\end{lem}

{\em Proof of Lemma}.   
Take  an arbitrary 
point $x \in Z_i$ and write $TV_x = TV_x' \oplus TV_x''$, where 
$\chi^{(i)}_1 \oplus \cdots \oplus \chi^{(i)}_{r_i}$ is the 
$H_i$-action on $TV_x'$ and $\tau^{(i)}_1  
\oplus \cdots \oplus \tau^{(i)}_{s_i}$
the $H_i$-action on 
$TV_x''$. By definition of $\overline{H_i}$, the group $\overline{H_i}$ 
acts trivially 
on $TV_x'$. Since ${\rm dim}\, TV_x' = r_i$ it is enough to show that 
$TV_x' = (TV_x)^{\overline{H_i}}$. Assume that $(TV_x)^{H_i'}$ is larger than 
$TV_x'$. Since $TV_x''$ is a $H_i$-invariant complement of 
$TV_x'$ in  $TV_x$, there exists a nonzero 
vector $v \in TV_x''$ such that $\overline{H_i}v = v$. 
Since $H_i/\overline{H_i}$ is abelian, 
the linear span of $Hv$ must split into $1$-dimensional 
irreducible representations of $H_i$. But all  
$\tau^{(i)}_1,  \ldots,   \tau^{(i)}_{s_i}$ have degree $\geq 2$. 
Contradiction.  \hfill $\Box$
\medskip

Let $W(V,G)$ be the union 
$W_1 \cup \cdots \cup W_k$. It is clear that $W(V,G)$ is $G$-invariant. 
We remark that $W_1, \ldots, W_k$ may not be necessary pairwise different.

{\sc Case 2a}:  $W(V,G)$ is smooth. 
Then we set 
$V_1$ to be the $G$-equivariant 
blow-up of $V_0$ 
with center $W(V,G)$. Using the same arguments as in Case 1, we conclude
that  $s(V_1, G) < s(V_0, G)$, because  the $H_i$-action in the normal 
space to $W_i$ at 
$x \in Z_i$ has no $1$-dimensional $H_i$-invariant subspaces. 
This allows to use the  induction on  $s(V, G)$.

{\sc Case 2b}:  $W(V,G)$ is not smooth. Then we must first 
resolve singularities
of $W(V,G)$. Denote by 
$k(x)$ the number of irreducible components of $W(V,G)$ containing $x \in 
W(V,G)$. 
Let $m(V,G) := \max_{x \in W(V,G)} k(x)$. We set   

\[ B(V,G) := \{ x \in W(V,G) \; :\;  k(x) = m(V,G) \} \subset 
W(V,G). \] 

It is clear that  $B(V,G)$ is $G$-invariant and smooth.  Now  we set 
$V_1$ to be the $G$-equivariant 
blow-up of $V_0$ with the center $B(V,G)$. 

We claim that $m(V_1,G) < m(V,G)$. This would follow immediately 
from the fact that 
all irreducible components of $W(V_1, G)$ are exactly the proper birational
transforms of the irreducible components of  $W(V, G)$.

Assume that  $W(V_1, G)$ contains an irreducible component $W_j'$
which is the connected component of $V_1^{\overline{H_j'}}$ containing a 
connected component $Z_j' \subset Z_{\rm max}(V_1,G)$, where 
$H_j' = St_G(x)$, $x \in Z_j'$. We 
necessary have $|H_j'| = s(V_1,G) = s(V,G)$. Then $\varphi_1(Z_j')$ 
must conside with a connected component $Z_i \subset Z_{\rm max}(V,G)$ 
and $H_j' = H_i$.
Let $E_i$ be the $\varphi_1$-exceptional divisor over $Z_i$,
$W_i$ the connected component  of  $V^{\overline{H_i'}}$
containing $Z_i$. Then  $B_i := \varphi_1(E_i)$ a connected component 
of $B(V,G)$. We have $Z_i
\subset B_i \subset W_i$ and $B_i \neq W_i$. It remains to observe that 
$Z_j'$ is contained in the proper birational transform 
 of $W_i$. The latter follows 
immediately from the decomposition of the linear representations of $H_i$ 
and $\overline{H_i}$ in the normal space to $B_i$ at 
any point $x \in Z_i$ (see the arguments in the proof  
of Lemma). 

Using the induction on  $m(V,G)$ we eventually come from Case 2b to Case2a.

\hfill $\Box$

\bigskip

I want to express my  gratitude  to Lev Borisov who explained  me an  
error in the abelianization algorithm presented in \cite{nonarch} \S 5.

\end{document}